\documentclass[11pt]{amsart}

\usepackage[applemac]{inputenc} 

\usepackage{amsthm,amssymb,amsbsy,amsmath,amsfonts,amssymb,amscd,mathrsfs}
\usepackage[normalem]{ulem}

\usepackage{mathtools} 
\mathtoolsset{showonlyrefs,showmanualtags} 

\long\def\symbolfootnote[#1]#2{\begingroup%
\def\thefootnote{\fnsymbol{footnote}}\footnote[#1]{#2}\endgroup} 

\usepackage{graphicx}
\DeclareGraphicsExtensions{.jpg}
\usepackage{color}
\usepackage{fullpage}
 \usepackage[usenames,dvipsnames]{pstricks}
 \usepackage{epsfig}
 \usepackage{pst-grad} 
 \usepackage{pst-plot} 
\usepackage[hypertex]{hyperref}
\newtheorem{theorem}{Theorem}
\newtheorem{corollary}[theorem]{Corollary}
\newtheorem{lemma}[theorem]{Lemma}
\newtheorem{proposition}[theorem]{Proposition}

\theoremstyle{definition} 
\newtheorem{definition}[theorem]{Definition}
\newtheorem{example}[theorem]{Example}
\theoremstyle{remark}
\newtheorem{remark}[theorem]{Remark}

\newcommand{\bt}{\begin{theorem}}
\newcommand{\et}{\end{theorem}}
\newcommand{\bl}{\begin{lemma}}
\newcommand{\el}{\end{lemma}}
\newcommand{\bp}{\begin{proposition}}
\newcommand{\ep}{\end{proposition}}
\newcommand{\bc}{\begin{corollary}}
\newcommand{\ec}{\end{corollary}}
\newcommand{\bdeff}{\begin{definition}}
\newcommand{\edeff}{\end{definition}}
\newcommand{\brem}{\begin{remark}}
\newcommand{\erem}{\end{remark}}
\newcommand{\bex}{\begin{example}}
\newcommand{\eex}{\end{example}}

\renewcommand{\r}[1]{(\ref{#1})}

\newcommand{\bR}{\mathbb{R}}
\newcommand{\be}{\begin{equation}}
\newcommand{\ee}{\end{equation}}

\newcommand{\lp}{\left(}
\newcommand{\rp}{\right)}
\newcommand{\lb}{\left[}
\newcommand{\rb}{\right]}

\newcommand{\Hess}{\nabla^{2}}

\DeclareMathOperator{\Cut}{Cut}
\DeclareMathOperator{\Con}{Con}

\renewcommand{\phi}{\varphi}

\newcommand{\bi}{\begin{itemize}}
\newcommand{\iii}{\item}
\newcommand{\ei}{\end{itemize}}
\newcommand{\bd}{\begin{description}}
\newcommand{\ed}{\end{description}}

\newcommand{\nn}{\nonumber}
\newcommand{\ba}[1]{\begin{array}{#1}}
\newcommand{\ea}{\end{array}}

\newcommand{\lam}{\lambda}

\newcommand{\g}{\gamma}
\newcommand{\al}{\alpha}
\newcommand{\eps}{\varepsilon}

\newcommand{\R}{\mathbb{R}}
\newcommand{\N}{\mathbb{N}}
\newcommand{\C}{\mathbb{C}}

\newcommand{\mb}[1]{\mathbb{ #1 }}
\newcommand{\mc}[1]{\mathcal{ #1 }}

\newcommand{\all}{\forall\,}

\newcommand{\la}{\langle}
\newcommand{\ra}{\rangle}

\newcommand{\virg}[1]{``#1''}
\newcommand{\tx}[1]{\mathrm{#1}}
\newcommand{\til}[1]{\widetilde{#1}}

\newcommand{\distr}{\mc{D}}

\newcommand{\metr}{\textsl{g}}

\newcommand{\Pg}[1]{\left\{ #1 \right\}}

\newcommand{\hp}{hypothesis}
\newcommand{\EXP}{\mathrm{exp}}
\newcommand{\Exp}{\mathrm{exp}}
\newcommand{\lapl}{\Delta}
\newcommand{\dive}{\text{div}}
\newcommand{\grad}{\nabla}

\newcommand{\HH}{\mc{H}}
\newcommand{\tcon}{t_{\mathrm{con}}}
\newcommand{\tcut}{t_{\mathrm{cut}}}

\newcommand{\smooth}{\phi}

\begin{document}
\begin{center} \noindent

\title{
Heat kernel asymptotics on sub-Riemannian manifolds with symmetries and applications to the bi-Heisenberg group}
\maketitle
Davide Barilari\\ 
{\footnotesize IMJ-PRG, Université Paris Diderot, 
 UMR CNRS 7586, Paris, France\\ {\tt davide.barilari@imj-prg.fr}}\\
\vskip 0.3cm
Ugo Boscain\\
{\footnotesize CNRS, CMAP Ecole Polytechnique, Paris, 
 Equipe INRIA GECO Saclay-\^Ile-de-France, Paris, France \\ {\tt ugo.boscain@polytechnique.edu}}
\vskip 0.3cm
Robert W. Neel
{\footnotesize \\Department of Mathematics, Lehigh University, Bethlehem, PA, USA\\ {\tt robert.neel@lehigh.edu}}

\vskip 0.3cm
\today
\end{center}

\vskip 0.2 cm

\begin{abstract}
By adapting a technique of Molchanov, we obtain the heat kernel asymptotics at the sub-Riemannian cut locus, when the cut points are reached by an $r$-dimensional parametric family of optimal geodesics. 
We apply these results to the bi-Heisenberg group, that is, a nilpotent left-invariant sub-Rieman\-nian structure on $\R^{5}$ depending on two real  parameters $\al_{1}$ and $\al_{2}$. We  develop some results about its geodesics and heat kernel associated to its sub-Laplacian and we illuminate some interesting geometric and analytic features appearing when one compares the isotropic ($\al_{1}=\al_{2}$) and the non-isotropic cases ($\al_{1}\neq \al_{2}$). In particular, we give the exact structure of the cut locus, and we get the complete small-time asymptotics for its heat kernel. \\

\end{abstract}

\section{Introduction}

The three-dimensional Heisenberg group is the simplest example of a sub-Riemannian manifold. It already possesses some features that are typical of more general sub-Riemannian structures: the Hausdorff dimension is bigger than the topological dimension, the cut locus and the conjugate locus starting from a point are adjacent to the point itself, spheres are not smooth, even for small radius, and the heat kernel small-time asymptotics on the diagonal are affected by the presence of the trivial abnormal minimizer (or, in other words, in sub-Riemannian geometry $q$ is not a smooth point for the exponential map based at $q$). 
The Heisenberg group is also the nilpotent approximation of any three-dimensional contact structure.

The natural next example, that shares with the Heisenberg group many of these properties, is the sub-Riemannian structure on $\R^{5}$ (where we consider coordinates $(x_{1},y_{1},x_{2},y_{2},z)$) 
where the following vector fields 
\begin{equation}
\begin{cases}
X_{1}=\partial_{x_1}-\dfrac{\al_1}{2}y_1\partial_z,  \qquad 
Y_{1}=\partial_{y_1}+\dfrac{\al_1}{2} x_1\partial_z, \\[0.3cm]
X_{2}=\partial_{x_2}-\dfrac{\al_2}{2}y_2\partial_z, \qquad
Y_{2}=\partial_{y_2}+\dfrac{\al_2}{2}x_2\partial_z.
\label{wectors}
\end{cases},\qquad \al_{1},\al_{2}\geq0, \ (\al_{1},\al_{2})\neq (0,0),
\end{equation}
define an orthonormal frame. Usually the case $\al_{1}=\al_{2}$ is called the five-dimensional Heisenberg group, while the case $\al_{1}\neq \al_{2}$ is referred as the \emph{non-isotropic} five-dimensional Heisenberg group, see \cite{bealsgaveaugreiner}. In the following we refer to the structure defined by \eqref{wectors} as the \emph{bi-Heisenberg group}.

When the bi-Heisenberg group is obtained as the nilpotent approximation of a sub-Riemannian structure with growth vector $(4,5)$, then the parameters  $\alpha_{1}$ and $\alpha_{2}$ keep track of some of the local structure of the original manifold. In other words, the presence of parameters is a consequence of the fact that the tangent space (or more precisely the nilpotent approximation) of a structure with growth vector $(4,5)$ depends on the point. (See \cite{corank1} for more details about the nilpotent approximation of a sub-Riemannian structure.)

Many features of the bi-Heisenberg group change when the two parameters become equal. The main purpose of this paper is to study these differences, in particular, as it concerns the structure of the cut locus and the behavior of the heat kernel.
An additional important case is the one in which one of the two parameters vanishes. This is the case in which the bi-Heisenberg structure is non-contact but is diffeomorphic to the product  $\mb{H}\times\R^2$, where $\mb{H}$ denotes the three-dimensional Heisenberg group. This case should be treated separately since abnormal minimizers are present.

The cut locus is obtained explicitly using known formulas for the cut time along each geodesics, see \cite{bealsgaveaugreiner,corank1}. The heat kernel asymptotics are obtained as a consequence of the main result of the paper, namely Theorem \ref{t-main}, which permits one to deduce the small-time heat kernel asymptotics at cut loci produced by symmetries of the exponential map. This result is obtained by refining a general method developed by Molchanov \cite{Molchanov} and the authors in \cite{srneel,imrn} (note in particular that it is a more detailed, sub-Riemannian version of Case 3.3 in \cite{Molchanov}). In the case of the bi-Heisenberg group, the cut points are also conjugate points and, as a consequence, one gets quite different asymptotic expansions in and out of the cut locus.

We note that an integral expression for the heat kernel of the bi-Heisenberg group is known, see \cite{bealsgaveaugreiner}. From this expression is not difficult to obtain an asymptotic expansion on the $z$-axis, but the extraction of the asymptotics on points other than these seems more difficult. We also note that stochastic approaches to small-time heat kernel asymptotics are possible, as in the recent work of  \cite{inahama}. While stochastic methods are quite powerful and provide insight into the entire diffusion associated to the sub-Laplacian, they are also technically demanding. In contrast, the Molchanov-inspired approach of the present paper treats conjugacy in a relatively simple, finite-dimensional way, and this makes it well suited to deriving explicit and geometrically meaningful asymptotics.

\subsection{Structure of the paper} Since the main result (Theorem \ref{t-main}) requires several preliminaries, we discuss first its application to the bi-Heisenberg group, in the next subsection.

After recalling basic definitions of sub-Riemannian geometry (Section \ref{s:bd}), in Section \ref{s:heat} we state and prove Theorem \ref{t-main}. Section \ref{s:2H} is devoted to the explicit computation of the optimal synthesis on the bi-Heisenberg group, together with the explicit formulas for the cut locus. Finally in Section \ref{s:final} we show that using the results of \cite{bealsgaveaugreiner} one can recover the expansion we found on a particular subset of the cut locus.

\subsection{Description of the results for the bi-Heisenberg group}

\bt[Structure of the cut locus]  \label{t-cut}The cut locus $\Cut_{0}$ starting from the origin for the bi-Heisenberg group (i.e., the sub-Riemannian metric on $\R^5$ for which  \r{wectors} is an orthonormal frame)  is characterized as follows:
\bi
\iii[(i)] if $0<\al_2=\al_1$ then $$\Cut_{0}=\{(0,0,0,0,z),~z\in\R\setminus\{0\} \},$$ 

\iii [(ii)]if $0<\al_2<\al_1$ then 
$$\Cut_{0}=\{(0,0,x_{2},y_{2},z),~|z|\geq  (x_2^2+y_2^2) K(\al_1,\al_2),~~(x_2,y_2,z)\in\R^3\setminus\{0\} \},$$ where
$$
K(\al_1,\al_2):= \frac{\al_{2}^{2}}{8}\left(\frac{2\pi}{\al_{1}}-\frac{1}{\al_{2}}\sin\left(2\pi\frac{\al_{2}}{\al_{1}}\right)\right)\sin^{-2}\left(\pi
\frac{\al_{2}}{\al_{1}}\right),
$$
\iii[(iii)] if $0=\al_2<\al_1$ then $$\Cut_{0}=\{(0,0,x_2,y_2,z),~(x_2,y_2,z)\in\R^3, z \neq 0 \}.$$
\ei

Moreover, in case (i)  every point of the cut locus is reached by a family of optimal geodesics parametrized by $S^{3}$; in cases (ii)-(iii)  every point of the cut locus is reached by a family of optimal geodesics parametrized by $S^{1}$. As a consequence all cut points are also conjugate points.

\et
Notice that the theorem above covers all the possible cases, up to exchanging the role of the indexes 1 and 2, which simply gives an isometric structure.

Notice also that in case (i) the cut locus is a one-dimensional manifold; in cases (ii)-(iii) the cut locus is a 3-dimensional manifold.
Moreover  $K(\al_1,\al_2)\to+\infty$ for $\al_2\to\al_1^{-}$ and  $K(\al_1,\al_2)\to 0$ for $\al_2\to0^{+}$. Hence  one can recover (i) and (iii) as limit cases of (ii). In all cases the cut locus is symmetric with respect to rotation along the $z$ axes.\\

In the following $d$ denotes the sub-Riemannian distance and $p_t(q_1,q_2)$ is the fundamental solution to the heat equation $\partial_t u=\lapl u$, where $\lapl$ is the sub-Laplacian defined on the bi-Heisenberg group (see Section \ref{s:bd} for precise definitions).

\bt[Asymptotics of the heat kernel]\label{t-heat}
 There exist positive constants $C_i$ (depending on $q$) such that for $t\to 0$,  the heat kernel for the bi-Heisenberg group satisfies
\bi
\iii[(i)] On-diagonal asymptotics:
\begin{equation}\nn
p_t(0,0)=\frac{C_1+O(t)}{t^3}.
\label{pt-diag}
\end{equation}
\iii[(ii)] If $q\notin \Cut_0\cup\{ 0\} $, then 
\begin{equation}\nn
p_t(0,q) = \frac{C_2+O(t)}{t^{5/2}}  e^{-d^{2}(0,q)/4t}.
\label{pt-nodiag-nocut}
 \end{equation}
\iii[(iii)] If $q\in \Cut_0$ and \\
\bi 
\iii[(a)]   $0<\al_2=\al_1$, then
\begin{equation}\nn
p_t(0,q) = \frac{C_3+O(t)}{t^4 } e^{-d^{2}(0,q)/4t},
\end{equation}
\iii[(b)]  $0\leq \al_2<\al_1$, then 
\label{pt-cut4}
\begin{equation}\nn
p_t(0,q) = \frac{C_4+O(t)}{t^3 } e^{-d^{2}(0,q)/4t}.
\label{pt-cut3}
\end{equation}
\ei 
\ei
\et 
Case (i) is well known and can be obtained from the explicit expression of the heat kernel or applying the result of Ben Arous \cite{benarouscut}. Notice that in this case one obtains $p_{t}(0,0)=t^{-Q/2}(C_{1}+O(t))$ where $Q=6$ is the Hausdorff dimension. Case (ii) is a consequence of a result of Ben Arous \cite{benarousdiag} and from the explicit expression of the cut locus (cf.\ Theorem \ref{t-cut}). Case (iii) is a consequence of Theorem \ref{t-main}. As explained in the introduction, the case in which $\al_2=0$ should be treated separately as a consequence of the presence of abnormal extremals.

\begin{remark} In Theorems \ref{t-cut} and \ref{t-heat} it is not restrictive to consider the origin as reference point, due to the left-invariance of the structure.  
\end{remark}

\begin{remark} With the same techniques, Theorems \ref{t-cut} and \ref{t-heat} can be generalized to the case of nilpotent structure of step 2 on $\R^{2\ell+1}$ with distribution of codimension 1. Notice that these structures depend on $\ell$ real parameters $\alpha_{1},\ldots,\alpha_{\ell}$. For an explicit description see \cite{corank1} (cf.\ also \cite{bealsgaveaugreiner}). In this case, the structure of the cut locus and the heat kernel asymptotics are determined by how many of these parameters coincide with the maximal one. 
\end{remark}

\section{Basic definitions}\label{s:bd}
In what follows $M$ is a connected orientable smooth manifold of dimension
$n\geq 3$.
\bdeff
A \emph{sub-Riemannian structure} on $M$ is a pair $(\distr,\metr)$,
where
\bi
\iii[(i)] $\distr$ is a smooth vector distribution of constant rank $k< n$
satisfying the \emph{H\"ormander condition}, i.e.\ a smooth map that
associates to $q\in M$  a  $k$-dimensional subspace $\distr_{q}$ of
$T_qM$ such that
\begin{equation} \label{Hor}
\qquad \text{span}\{[X_1,[\ldots[X_{j-1},X_j]]]_{q}~|~X_i\in\mathfrak{X}(\distr),\,
j\in \N\}=T_qM, \ \all q\in M,
\end{equation}
where $\mathfrak{X}(\distr)$ denotes the set of \emph{horizontal smooth
vector fields} on $M$, i.e.\
$$\mathfrak{X}(\distr)=\Pg{X\in\mathrm{Vec}(M)\ |\ X(q)\in\distr_{q}~\
\forall~q\in M}.$$
\iii[(ii)] $\metr_q$ is a Riemannian metric on $\distr_{q}$ which is smooth
as function of $q$. We denote  the norm of a vector $v\in \distr_{q}$
by 
$|v|_{\metr}=\sqrt{\metr_{q}(v,v)}.$
\ei
\edeff
A Lipschitz continuous curve $\g:[0,T]\to M$ is said to be
\emph{horizontal} (or \emph{admissible}) if
$$\dot\g(t)\in\distr_{\g(t)}\qquad \text{ for a.e.\ } t\in[0,T].$$
The {\it length} of  an horizontal curve $\g:[0,T]\to M$ is
\begin{equation}
\label{e-lunghezza}
\ell(\g)=\int_0^T |\dot{\g}(t)|_{\metr}~dt.
\end{equation}
Notice that $\ell(\g)$ is invariant under time reparametrization of $\g$.
The {\it sub-Riemannian (or Carnot-Carath\'eodory) distance} on $M$ is the
function
\begin{equation}
\label{e-dipoi}
d(q_0,q_1)=\inf \{\ell(\g)\mid \g(0)=q_0,\g(T)=q_1, \g\ \mathrm{horizontal}\}.
\end{equation}
The \hp\ of connectedness of $M$ and the H\"ormander condition
guarantee the finiteness and the continuity of $d:M\times M\to \R$ with
respect to the topology of $M$ (Chow-Rashevsky theorem, see for
instance \cite{agrachevbook}). The function $d$ defined in \eqref{e-dipoi} gives to $M$ the structure
of a metric space (see \cite{agrachevbook}) compatible with the original topology of $M$. The sub-Riemannian manifold is said {\it complete} if $(M,d)$ is complete as a metric space.

Locally, the pair $(\distr,\metr)$ can be specified by assigning a set of $k$ smooth vector fields spanning $\distr$, being an orthonormal frame for $\metr$, i.e.  
\begin{equation}
\label{trivializable}
\distr_{q}=\text{span}\{X_1(q),\dots,X_k(q)\}, \qquad \qquad \metr_q(X_i(q),X_j(q))=\delta_{ij}.
\end{equation}
In this case, the set $\Pg{X_1,\ldots,X_k}$ is called a \emph{local orthonormal frame} for the sub-Riemannian metric. 

When a local orthonormal frame is defined globally, the sub-Riemannian manifold is called {\it trivializable}.

The sub-Riemannian metric can also be expressed locally in ``control form'' as follows. Fix $\Pg{X_1,\ldots,X_k}$  a {local orthonormal frame} and consider the control system:
\begin{equation}
\dot q=\sum_{i=1}^k u_i X_i(q).
\end{equation}
The problem of finding the shortest curve minimizing that joins two fixed points $q_0,~q_1\in M$ is naturally formulated as the optimal control problem of minimizing
\begin{equation} \label{eq:lnonce}
~~~\int_0^T  \sqrt{
\sum_{i=1}^k u_i^2(t)}~dt\, \quad \text{subject to} \quad q(0)=q_0,~~~q(T)=q_1.
\end{equation}

A sub-Riemannian structure on $M$ is said to be \emph{left-invariant} if $M=G$, where $G$ is a Lie group, and both $\distr$ and $\metr$ are left-invariant over $G$.
Left-invariant sub-Riemannian manifolds are trivializable. An example of left-invariant sub-Riemannian manifold is provided by the bi-Heisenberg group presented in the introduction (cf.\ also Section \ref{s:2H}).

Now we briefly recall some facts about sub-Riemannian
geodesics. In particular, we define the
sub-Riemannian exponential map.

\bdeff A \emph{geodesic} for a sub-Riemannian manifold
$(M,\distr,\metr)$ is an admissible curve $\g:[0,T]\to M$ such that
 $|\dot{\g}(t)|_{\metr}$  is constant and, for every $t\in [0,T]$ there exists $\eps>0$ such that the restriction
$\g_{|_{[t-\eps,t+\eps]}}$ is a minimizer of $\ell(\cdot)$.
A geodesic satisfying $|\dot{\g}(t)|_{\metr}=1$ is said to be
parameterized by arclength.

\edeff

A version of the  Hopf-Rinow theorem (see \cite[Chapter 2]{burago})
 ensures that  if the sub-Riemannian manifold is complete as metric space, then all geodesics are defined for every $t\geq0$ and that for every two points there exists a minimizing geodesic
connecting them.

Left invariant sub-Riemannian manifolds are always complete.

Trajectories minimizing the distance between two points are solutions
of first-order necessary conditions for optimality, which in the case
of sub-Riemannian geometry are described as follows.
\bt\label{t:pmpw}
Let $q(\cdot):t\in[0,T]\mapsto q(t)\in M$ be a solution of the
minimization problem \eqref{eq:lnonce} such that
$|\dot q(t)|_{\metr}$  is constant, and let $u(\cdot)$ be the
corresponding control.
Then there exists a Lipschitz map $p(\cdot): t\in [0,T] \mapsto
p(t)\in T^{*}_{q(t)}M\setminus\{0\}$  such that one and only one of
the following conditions holds:
\bi
\iii[(i)]
$
\dot{q}=\dfrac{\partial H}{\partial p}, \quad
\dot{p}=-\dfrac{\partial H}{\partial q}, \quad
u_{i}(t)=\la p(t), X_{i}(q(t))\ra,
\\$
where $H(q,p)=\frac{1}{2} \sum_{i=1}^{k} \la p,X_{i}(q)\ra^{2}$.
\vspace{0.2cm}
\iii[(ii)]
$
\dot{q}=\dfrac{\partial \HH}{\partial p}, \quad
\dot{p}=-\dfrac{\partial \HH}{\partial q}, \quad
0=\la p(t), X_{i}(q(t))\ra,
\\$
where $\HH(t,q,p)=\sum_{i=1}^{k}u_{i}(t) \la p,X_{i}(q)\ra$.
\ei
\et
\noindent
For an elementary proof of Theorem \ref{t:pmpw} see \cite{nostrolibro}.

\brem If $(q(\cdot),p(\cdot))$ is a solution of (i) (resp.\ (ii)) then
it is called a \emph{normal extremal} (resp.\ \emph{abnormal
extremal}). It is well known that if $(q(\cdot),p(\cdot))$ is a normal
extremal then $q(\cdot)$ is a geodesic (see
\cite{nostrolibro,agrachevbook}). This does not hold in general for
abnormal extremals. An admissible trajectory $q(\cdot)$ can be at the
same time normal and abnormal (corresponding to different covectors).
If an admissible trajectory $q(\cdot)$ is normal but not abnormal, we
say that it is \emph{strictly normal}.

\erem

\bdeff
A minimizer $\gamma:[0,T]\to M$ is said to be \emph{strongly normal}
if for every $[t_1,t_2]\subset [0,T]$, $\gamma|_{[t_1,t_2]}$ is not an
abnormal minimizer.\edeff

In the following we denote by $(q(t),p(t))=e^{t\vec{H}}(q_{0},p_{0})$
the solution of $(i)$ with initial condition
$(q(0),p(0))=(q_{0},p_{0})$. Moreover we denote by $\pi:T^{*}M\to M$
the canonical projection.

 Normal extremals (starting from $q_{0}$) parametrized by arclength
correspond to initial covectors $p_{0}\in \Lambda_{q_{0}}:=\{p_{0}\in
T^{*}_{q_{0}}M | \, H(q_{0},p_{0})=1/2\}.$
\bdeff  \label{d:exp}
Let $(M,\distr,\metr)$ be a complete sub-Riemannian manifold and
$q_0\in M$. We define the \emph{exponential map} starting from $q_{0}$
as
\begin{equation} \label{eq:expmap}
\EXP_{q_{0}}: \Lambda_{q_{0}}\times \R^{+} \to M, \qquad
\EXP_{q_{0}}(p_{0},t)= \pi(e^{t\vec{H}}(q_{0},p_{0})).
\end{equation}
\edeff

Next, we recall the definition of cut and conjugate time.

\bdeff \label{def:cut} Let $q_{0}\in M$ and
$\g(t)$ an arclength geodesic starting from $q_{0}$.
The \emph{cut time} for $\g$ is $\tcut(\g)=\sup\{t>0,\, \g|_{[0,t]}
\text{ is optimal}\}$. The \emph{cut locus} from $q_{0}$ is the set
$\Cut_{q_{0}}=\{\g(\tcut(\g)), \g$  arclength geodesic from
$q_{0}\}$.
\edeff

\bdeff \label{def:con} Let $q_{0}\in M$ and
$\g(t)$ a normal arclength geodesic starting from $q_{0}$ with initial
covector $p_{0}$. Assume that $\g$ is not abnormal.
\begin{itemize}
\item[(i)] The \emph{first conjugate time} of $\g$ is
$\tcon(\g)=\min\{t>0,\  (p_{0},t)$ is a critical point of
$\EXP_{q_{0}}\}$. 
\item[(ii)] The \emph{first conjugate locus} from $q_{0}$ is the
set $$\Con_{q_{0}}=\{\g(\tcon(\g)), \g \  \text{arclength geodesic from}\ 
q_{0}\}.$$
\end{itemize}
\edeff
It is well known that, for a geodesic $\g$ which is not abnormal, the
cut time $t_{*}=\tcut(\g)$ is either equal to the conjugate time, or
there exists another geodesic $\til{\g}$ such that
$\g(t_{*})=\til{\g}(t_{*})$ (see for instance \cite{nostrolibro,agrexp}).

\brem
In sub-Riemannian geometry, the  map $\Exp_{q_{0}}$ is never a local
diffeomorphism at $0\in T^{*}_{q_{0}}M$. As a consequence one can show that the
sub-Riemannian balls are never smooth and both the cut and the
conjugate loci from $q_{0}$ are adjacent to the point $q_{0}$ itself
(see \cite{agratorino,nostrolibro}).

\erem
\subsection{The sub-Laplacian}\label{s:lapl}
In this section we define the sub-Riemannian Laplacian
 on a sub-Riemannian manifold $(M, \distr, \metr)$ endowed with a
smooth volume $\mu$.

The sub-Laplacian is the natural generalization of the
Laplace-Beltrami operator defined on a Riemannian manifold, defined as
the divergence of the gradient.
The \emph{horizontal gradient} is the
unique operator $\grad: \mathcal{C}^\infty(M)\to \mathfrak{X}(\distr)$ defined by $$
\metr_q(\grad\phi(q),v)=d\phi_{q}(v), \quad \all \phi \in
\mc{C}^{\infty}(M), \, q\in M,~v\in \distr_q.$$
By construction, it is a horizontal vector field. If
$X_{1},\ldots,X_{k}$ is a local orthonormal frame, it is easy to see
that it is written as follows
$\grad \phi= \sum_{i=1}^{k} X_i(\phi) X_i,$
where $X_i(\phi)$ denotes the Lie derivative of $\phi$ in the
direction of $X_i$.

The divergence of a vector field $X$ with respect to a volume $\mu$ is
the function $\dive_{\mu}X$ defined by the identity $L_{X}\mu=(\dive_{\mu} X
)\mu$, where $L_{X}$ stands for the Lie derivative with respect to $X$.

The sub-Laplacian associated with  the sub-Riemannian structure, i.e.\
$\lapl \phi=\dive_{\mu}(\grad \phi),$
 is written in a local orthonormal frame $X_{1},\ldots,X_{k}$ as follows
\begin{equation} \label{eq:lapldiv}
\lapl= \sum_{i=1}^{k}X_{i}^{2}+ (\dive_{\mu} X_{i}) X_{i}.
\end{equation}
Notice that $\lapl$ is always expressed as the sum of squares of the
element of the orthonormal frame plus a first order term that belongs
to the distribution and depends on the choice of the volume $\mu$.

In what follows we denote  $p_t$  the heat kernel associated to the heat equation $\partial_t u=\lapl u$. Its existence, smoothness, symmetry, and positivity are guaranteed by classical results since the operator \eqref{eq:lapldiv} is in divergence form and $M$ is complete, see for instance \cite{strichartz}.

\brem \label{r:group}
In the case of a left-invariant structure on a Lie group (and in
particular for a nilpotent structure), an intrinsic choice of volume is provided by the left Haar measure. When the Lie group is unimodular, which is the case in particular for nilpotent groups, the left Haar measure is both left and right invariant. Left invariant vector fields on a Lie group $G$ are divergence free with respect to the right Haar measure. As a consequence the sub-Laplacian associated to the Haar measure has the form of \virg{sum of squares} 
(see also \cite{hypoelliptic})
$$\lapl=\sum_{i=1}^k X_i^2.$$
For unimodular left invariant structures one can show that the Haar measure coincides with the Popp volume, that is a well defined measure for equiregular sub-Riemannian manifolds (see \cite{BRpopp}).
\erem

\section{Small time asymptotic of the heat kernel}\label{s:heat}

We begin by recalling the general framework for computing asymptotics developed in \cite{srneel}, which is built upon the basic insight of Molchanov \cite{Molchanov}. Let $M$ be a complete sub-Riemannian (or Riemannian) manifold of dimension $n$. Consider distinct points $x$ and $y$ such that every minimizing geodesic connecting them is strongly normal. (In this section, we will use $x$, $y$ and $z$ to denote points of $M$, to avoid the need for subscripts as in $q_0$, $q_1$, etc. Because we work on a general sub-Riemannian manifold, there should be no confusion with using these letters to denote Euclidean coordinates on $\bR^5$, as when we specifically discuss the bi-Heisenberg group.) Let $\Gamma\subset M$ be the set of midpoints of minimizing geodesics from $x$ to $y$ (or from $y$ to $x$, the situation is symmetric), and let $N(\Gamma)$ denote a neighborhood of $\Gamma$. Further, let $h_{x,y}(z)=\lb d^2(x,z)+d^2(z,y)\rb/2$; then $h_{x,y}$ attains its minimum of $d^2(x,y)/4$ exactly on $\Gamma$ (see Lemma 21 of \cite{srneel}). Because $h_{x,y}$ is continuous and sub-Riemannian balls are compact (by completeness of $M$), it follows that $\Gamma$ is compact. Further, since strong normality is an open condition and $\Gamma$ is a positive distance from $x$, $y$, and their cut loci, $h_{x,y}$ is smooth in a neighborhood of $\Gamma$. We assume that $N(\Gamma)$ is small enough so that its closure is positive distance from $x$, $y$, and their cut loci and that $h_{x,y}$ is smooth on its closure.

Continuing, suppose we give $M$ a sub-Laplacian $\lapl$ coming from a smooth volume $\mu$ as in Section \ref{s:lapl}.
Let $\Sigma\subset M\times M$ be the set of $(x,y)$ such that $x\neq y$ and there exists a unique, strongly normal, non-conjugate minimizing geodesic connecting $x$ and $y$. Then from the results contained in Ben Arous \cite{benarouscut}, it follows that for $(x,y)\in\Sigma$,
\[
p_t(x,y) = \frac{1}{t^{n/2}}\exp\lp -\frac{d^2(x,y)}{4t}\rp\lp c_0(x,y)+ O(t)\rp \quad\text{as $t\rightarrow0$},
\]
where $c_0(x,y)$ is smooth and positive on $\Sigma$ and the $O(t)$ is uniform over compact subsets of $\Sigma$.

We can now give the basic, general  result for determining small-time asymptotics at the cut locus. The first part of Theorem 27 of \cite{srneel} gives the following.

\bt\label{THM:GenAsymptotics}
Let $M$ be an $n$-dimensional (complete) sub-Riemannian manifold with smooth measure $\mu$ and corresponding heat kernel $p_t$, and let $x$ and $y$ be distinct points such that all minimal geodesics from $x$ to $y$ are strongly normal. Then with the above notation, for any sufficiently small $N(\Gamma)$, we have
\[
p_t(x,y)= \int_{N(\Gamma)}\lp\frac{2}{t}\rp^n e^{-h_{x,y}(z)/t}\lp c_0(x,z)c_0(z,y)+ O(t)\rp\mu(dz) ,
\]
where the $O(t)$ term is uniform over $N(\Gamma)$.
\et

We wish to apply this theorem in concrete cases where we have an explicit expression for $h_{x,y}$, and possibly for $c_0$. For what follows, we will need the case when $\Gamma$ is an embedded (necessarily compact) $r$-dimensional smooth submanifold of $M$ and the Hessian of $h_{x,y}$ restricted to normal bundle of $T\Gamma$ in $TM$ is non-degenerate (because in general there is no Riemannian metric on $M$, there is no notion of a canonical representative of the normal bundle as an orthogonal complement, but non-degeneracy of the Hessian is nonetheless well defined). In other words, $h_{x,y}$ is a Morse-Bott function with critical set $\Gamma$. Note that because the Hessian of $h_{x,y}$ at $z\in\Gamma$ is always non-degenerate along the minimal geodesic from $x$ to $y$ through $z$, we must have $0\leq r<n$.

In this situation, for any local coordinates $(u_1,\ldots,u_r)$ on an open subset $U$ of $\Gamma$, we can extend them to local coordinates $(u_1,\ldots,u_r,u_{r+1},\ldots,u_n)$ on an open subset $V$ of $M$ such that $\Gamma$ is locally given by $u_{r+1}=\cdots=u_n=0$ and $U=V\cap \Gamma$. Further, perhaps by shrinking $U$ and $V$, we can assume that these coordinates can be extended to a neighborhood of $\overline{V}$. We call such a local coordinate system on $M$ adapted to $\Gamma$. Then $\mu$ has a smooth density $F(u)$ with respect to $du_1\cdots du_n$, and $F(u)$ along all with of its derivatives in the $u_i$ are bounded on $V$. Next, note that for $z\in U$, $\partial_{u_i}\partial_{u_j} h_{x,y} (z)=0$ if either $i\leq r$ or $j\leq r$. Then let $H^c h_{x,y}(z)$ be the $(n-r)\times(n-r)$ matrix with $i,j$ entry $\partial_{u_{r+i}}\partial_{u_{r+j}} h_{x,y} (z)$ (this depends on the coordinate system $u$, but we supress that in the notation). Note that all the entries of $H^c h_{x,y}(z)$ along with all of their derivatives in the $u_i$ are bounded on $U$. Now we claim that there is measure $\nu_{\Gamma}^{\mu,h_{x,y}}$ on $\Gamma$, induced by $\mu$ and the Morse-Bott function $h_{x,y}$, that is defined by having the following local expression in any local adapted coordinates,
\[
\nu_{\Gamma}^{\mu,h_{x,y}} = \frac{F(u_1,\ldots,u_r,0,\ldots,0)}{\sqrt{\det H^c h_{x,y}(u_1,\ldots,u_r)}} du_1\ldots du_r.
\]
 
To see that this is a globally well-defined measure on $\Gamma$, first note that this expression depends only on the local coordinates $(u_1,\ldots,u_r)$ on $\Gamma$ and not on the extension to adapted coordinates on $M$; more precisely, the factor of $F/\sqrt{\det H^c h_{x,y}}$ is invariant under change of fiber coordinates on the normal bundle. This is because if we change the coordinates $(u_{r+1},\ldots,u_n)$ on any fixed fiber of the normal bundle to some other adapted choice, both $F^{\Gamma}$ and $\sqrt{\det H^c h_{x,y}}$ are multiplied by the aboslute value of the determinant of the Jacobian, which clearly cancels. Further, if we change the coordinates $(u_1,\ldots,u_r)$ on $\Gamma$, it is clear that the given expression for $\nu_{\Gamma}^{\mu,h_{x,y}}$ transforms like a measure, and this is enough to establish the claim.

The motivation for introducing $\nu_{\Gamma}^{\mu,h_{x,y}}$ is that it encodes the relevant behavior of both $\mu$ and the Hessian of $h_{x,y}$ needed for the following asymptotic expansion.

\bt\label{t-main}
Let $M$ be an $n$-dimensional complete sub-Riemannian manifold provided
with a smooth volume $\mu$ and let $p_{t}$ be the heat kernel of the
sub-Riemannian heat equation. Let $x$ and $y$ be distinct and assume that every optimal geodesic joining $x$
to $y$ is strongly normal. Define
$$
\mathcal{O}:=\{\bar p\in \Lambda_x  ~|~ \EXP_{x}(\bar p,d(x,y))=y\}
$$
Assume that:
\bi
\item[(i)] $\mathcal{O}$ is a submanifold of $\Lambda_x$ of dimension $r$. 
\item[(ii)]  for every $\bar p\in \mathcal{O}$ we have 
$
\dim \ker 
D_{\bar p,d(x,y) }  \EXP_{x}=r.
$
\ei
Then there exists a positive
constant $C$ such that
\begin{equation}
p_t(x,y) = \frac{C+O(t)}{t^{\frac{n+r}{2}}}  e^{-d^{2}(x,y)/4t}, \qquad \text{
as } t\rightarrow 0.
\end{equation}
Further, $\EXP_{x}(\cdot,d(x,y)/2)$ maps $\mathcal{O}$ diffeomorphically onto $\Gamma$, and the Hessian of $h_{x,y}$ restricted to the normal bundle of $T\Gamma$ in $TM$ is non-degenerate, so that $\nu_{\Gamma}^{\mu,h_{x,y}}$ is defined as above. Finally, the constant $C$ in the expansion is given by
\begin{equation}
C = 2^{(3n-r)/2}\pi^{(n-r)/2} \int_{z\in \Gamma} c_0(x,z)c_0(z,y) \nu_{\Gamma}^{\mu,h_{x,y}} \lp dz\rp  .
\end{equation}
\et

\begin{proof} By its definition, $\Gamma$ is the image of  $\mathcal{O}$ under $\EXP_x( \cdot, d(x,y)/2)$. Take $\bar p \in \mathcal{O}$. By the assumption of strong normality and basic properties of the exponential map, $\EXP_x( \cdot, d(x,y)/2)$ is a diffeomorphism from a neighborhood of $\bar p$ in $\Lambda_x$ to a neighborhood of $z_0=\EXP_x( \bar p, d(x,y)/2)$ in $M$. Thus $\Gamma$ is locally (near $z_0$) an immersed submanifold of $M$ of dimension $r$.  Further, because the geodesics corresponding to covectors in $\Lambda_x$ are optimal until time $d(x,y)$ and the cut locus of $x$ is closed, $\EXP_x( \bar p, d(x,y)/2)$ is injective on a neighborhood of $\mathcal{O}$. It follows that $\Gamma$ is a compact embedded submanifold of $M$ of dimension $r$, that $\EXP_x( \bar p, d(x,y)/2)$ is a diffeomorphism from a neighborhood of $\mathcal{O}$ in $\Lambda_x$ to a neighborhood of $\Gamma$ in $M$ that restricts to a diffeomorphism from $\mathcal{O}$ to $\Gamma$, and that $\mathcal{O}$ is compact.

Since $h_{x,y}$ is constant on $\Gamma$ (see also Theorem 24 of \cite{srneel}) and $h_{x,y}$ is smooth near $\Gamma$, the tangent space to $\Gamma$ lies in the kernel of $\Hess h_{x,y}$. Since this tangent space has dimension $r$, we see that this tangent space is in fact equal to the kernel of $\Hess h_{x,y}$ (at every point of $\Gamma$). In particular, $\Hess h_{x,y}$ is non-degenerate on the $(n-r)$-dimensional normal space to $\Gamma$. So we are in the situation discussed above, and we can define $\nu_{\Gamma}^{\mu,h_{x,y}}$.

Continuing, by general facts about compact submanifolds (like the tubular neighborhood theorem), it follows that we can find a sufficiently small $N(\Gamma)$ so that Theorem \ref{THM:GenAsymptotics} holds for it and such that $N(\Gamma)$ is the union of a finite number of open sets $V_1,\ldots, V_k$ with the following properties:
\begin{itemize}
\item There are coordinates $(u^i_1,\ldots, u^i_n)$ on each $V_i$ which form a set of coordinate charts for $N(\Gamma)$, and these coordinates can be extended to a neighborhood of $\overline{V_i}$.
\item  If $U_i=V_i\cap\Gamma$, then each $U_i$ is nonempty and corresponds to the subset $\{u^i_{r+1}=\cdots=u^i_n=0\}$ of $V_i$. In particular, $(u^i_1,\ldots, u^i_n)$ are local coordinates that are adapted to $\Gamma$, as above.
\item Each $V_i$ can be written, in the coordinates $(u^i_1,\ldots, u^i_n)$, as the product $U\times(-\eps,\eps)^{n-r}$ for some $\eps>0$.
\item There exists a partition of unity $\phi_1,\ldots,\phi_k$ for $\Gamma$ subordinate to the covering $U_1,\ldots, U_k$.
\end{itemize}

In the computations that follow, we drop the superscript $i$ from the coordinates $(u^i_1,\ldots, u^i_n)=(u_1,\ldots, u_n)$ on each $V_i$, in order to lighten the notation. We can extend each $\phi$ to be a function on $V_i$ by letting them be constant on the normal fibers (which are given by fixing $u_1,\ldots,u_r$); then $\sum_{i=1}^k \phi_i =1$ on all of $N(\Gamma)$. Now we can apply Theorem \ref{THM:GenAsymptotics} and write the integral over $N(\Gamma)$ in terms of this system of coordinates charts to find
\begin{equation}\label{Eqn:BasicAsyIntegral}\begin{split}
p_t(x,y)= &\lp\frac{2}{t}\rp^n \sum_{i=1}^k \int_{U_i}\Bigg\{ \phi_i(u_1,\ldots,u_r)  \\
& \times \int_{(-\eps,\eps)^{n-r}} 
 e^{-h_{x,y}(u_1,\ldots,u_n)/t} \big(c_o(x,u_1,\ldots,u_n)c_0(u_1,\ldots,u_n,y) + O(t)\big)  \\
&  \times F_i(u_1,\ldots,u_n)
du_{r+1}\cdots du_n \Bigg\} du_1\cdots du_r
\end{split}\end{equation}
where the $O(t)$ terms are uniform and $F_i$ is the density of $\mu$ with respect to $du_1\cdots du_n$.
The asymptotic behavior of the inner integrals can be computed from the standard Laplace asymptotic formula (see page 198 of \cite{KanwalEstrada}), since the Hessian of $h_{x,y}$ in the normal directions is not degenerate, and this behavior is uniform by the smoothness and boundedness of everything involved. This gives 
\[\begin{split}
p_t(x,y)&= \lp\frac{2}{t}\rp^n \sum_{i=1}^k \int_{U_i} \phi_i(u_1,\ldots,u_r)
e^{-d^2(x,y)/4t}    \lp2\pi t\rp^{(n-r)/2} \\
&\times (c_0|_{\Gamma}(x,u_1,\ldots,u_r)c_0|_{\Gamma}(u_1,\ldots,u_r,y)+O(t)) \nu_{\Gamma}^{\mu,h_{x,y}} \lp du_1\cdots du_r \rp .
\end{split}\]
Using that the $\phi_i$ are a partition of unity and that the $O(t)$ term is uniform and thus its integral is also $O(t)$, we re-write this as
\[
p_t(x,y) = \frac{2^{(3n-r)/2}\pi^{(n-r)/2}}{t^{(n+r)/2}} e^{-d^2(x,y)/4t} \lb 
\int_{z\in \Gamma}c_0(x,z)c_0(z,y)\nu_{\Gamma}^{\mu,h_{x,y}} \lp dz\rp +O(t)\rb .
\]
It only remains to see that the integral over $\Gamma$ is positive (so that the $C$ in the theorem is positive), but this follows immediately from the positivity of $c_0$ and the fact that $\nu_{\Gamma}^{\mu,h_{x,y}}$ clearly gives positive measure to every non-empty open (Borel) subset of $\Gamma$. \end{proof}

\begin{remark}
The statement of the theorem makes it clear that if one is only interested in the power of $t$ appearing in the asymptotics of $p_t(x,y)$, which is the quantity that is most related to the (conjugacy) structure of minimal geodesics from $x$ to $y$, then only the properties $(i)$ and $(ii)$ of the exponential map are needed. It is only if one wants to determine the constant $C$ that any specific control of $\mu$, $c_0(\cdot,\cdot)$, or the Hessian of $h_{x,y}$ on the normal bundle are needed.
\end{remark}

\begin{remark}
This theorem covers the most natural case, in which the dimension of $\mathcal{O}$ is equal to the number of independent directions in which the optimal geodesics are conjugate, i.e.\ to $\dim \ker D_{\bar p,d(x,y)}\EXP_{x}$. Not only does this seem most common in situations of interest, but also there is a clear relationship between the geodesic structure and the behavior of $h_{x,y}$. In particular, the non-conjugacy in the normal directions, to $\mathcal{O}$ or equivalently $\Gamma$ (in what follows we freely use the fact that the exponential map gives a diffeomorphism between a neighborhood of $\mathcal{O}$ and a neighborhood of $\Gamma$), corresponds exactly to the non-degeneracy of the Hessian of $h_{x,y}$ in those directions. More concretely, the non-degeneracy of the Hessian along the normal fibers means that we can find a further change of coordinates $(u_{r+1},\ldots,u_n)\mapsto (\tilde{u}_{r+1},\ldots,\tilde{u}_n)$ on each fiber such that $h_{x,y}$, restricted to the fiber, is given by the sum of squares $\tilde{u}^2_{r+1}+\ldots+\tilde{u}^2_n$. Indeed, this normal form for $h_{x,y}$ is what underlies the asymptotic expansion of the ``inner integrals'' in \eqref{Eqn:BasicAsyIntegral}.

This indicates that the above computation of the asymptotic expansion can be extended to other normal forms for $h_{x,y}$ along the normal fibers. To briefly illustrate this, suppose that we drop condition (ii) in Theorem \ref{t-main} and instead assume that near each $z\in\Gamma$ we can find local coordinates such that $h_{x,y}= u^2_{r+1}+\ldots+u^2_{n-1}+u^4_n$. This implies that $\dim \ker D_{\bar p,d(x,y)}\EXP_{x}=r+1$ for $\bar p\in\mathcal{O}$. Moreover, as discussed in \cite{imrn} (see especially Lemmas 27 and 29), if the exponential map restricted to a normal fiber has a singularity of type $A_3$ in the Arnol'd classification, then $h_{x,y}$ will have this normal form when restricted to the same normal fiber. So there is again a connection between the geodesic geometry and the behavior of $h_{x,y}$, even if it is not as clean. Nonetheless, once we assume this normal form for $h_{x,y}$, the argument proceeds in precisely the same way, except that the ``inner integrals'' in \eqref{Eqn:BasicAsyIntegral} now have asymptotics determined by this different local expression for $h_{x,y}$. Such an expansion is again well known, and completing the computation, one finds that 
\[
p_t(x,y) = \frac{C+O\lp \sqrt{t}\rp}{t^{\frac{n+r+1}{2}-\frac{1}{4}}}  e^{-d^{2}(x,y)/4t}, \qquad \text{
as } t\rightarrow 0.
\]
Here $C$ is again some positive constant, the precise value of which is given by integrating some local data (depending on $c_0$, the volume measure and behavior of $h_{x,y}$) over $\Gamma$.

\end{remark}

\section{The bi-Heisenberg group} \label{s:2H}

Recall that the bi-Heisenberg group is the sub-Riemannian structure on $\R^{5}$ (where we consider coordinates $(x_{1},y_{1},x_{2},y_{2},z)$) 
where the following vector fields 
\begin{equation}
\begin{cases}
X_{1}=\partial_{x_1}-\dfrac{\al_1}{2}y_1\partial_z,  \qquad 
Y_{1}=\partial_{y_1}+\dfrac{\al_1}{2} x_1\partial_z, \\[0.3cm]
X_{2}=\partial_{x_2}-\dfrac{\al_2}{2}y_2\partial_z, \qquad
Y_{2}=\partial_{y_2}+\dfrac{\al_2}{2}x_2\partial_z.
\label{wectors2}
\end{cases},\qquad \al_{1},\al_{2}\geq0, \ (\al_{1},\al_{2})\neq (0,0),
\end{equation}
define an orthonormal frame. If we set $Z=\partial_{z}$, it is easy to see that the only non trivial commutation relations are
$$[X_{1},Y_{1}]=\al_{1}Z,\qquad [X_{2},Y_{2}]=\al_{2}Z,$$ hence the structure is bracket generating.
Notice that it is a nilpotent left-invariant structure (indeed it is a Carnot group). Moreover the structure is contact  if and only if $\al_{1},\al_{2}>0$. When one of the two parameters is zero the sub-Riemannian structure is the product  $\R^{2}\times \mathbb{H}$ where $\mathbb{H}$ is the 3-dimensional Heisenberg group. \\

\subsection{Exponential map and synthesis}
In this section we compute explicitly the cut locus starting from the origin in the bi-Heisenberg group. In this section we assume $0\leq \al_{2}\leq\al_{1}$, treating the case $\al_{2}=0$ separately. The case $0\leq \al_{1}\leq\al_{2}$ can be obtained by exchanging the role of the indexes 1 and 2.


\subsubsection{Contact case} 

In the contact case (i.e., $\al_{2}>0$) there are no non-constant abnormal extremals. One then reduces the computations of the extremals (see Theorem \ref{t:pmpw}) to the solution of the Hamiltonian system defined by the Hamiltonian
$$H(p,q)=\frac{1}{2}\sum_{i=1}^{2}\la p,X_{i}(q)\ra^{2}+\la p,Y_{i}(q)\ra^{2}.$$
The arclength geodesics starting from the origin are parametrized by the initial covector $p_{0}\in T^{*}_{0}M$ belonging to the level set $\Lambda_{0}:=H^{-1}(1/2)\cap T^{*}_{0}M$. If $(p_{x_{1}},p_{y_{1}},p_{x_{2}},p_{y_{2}},w)$ are the dual variables to $(x_{1},y_{1},x_{2},y_{2},z)$ then 
$$\Lambda_{0}=\{(p_{x_{1}},p_{y_{1}},p_{x_{2}},p_{y_{2}},w)\,|\, p^{2}_{x_{1}}+p^{2}_{y_{1}}+p^{2}_{x_{2}}+p^{2}_{y_{2}}=1\}\simeq S^{3}\times \R$$
Performing the change of variable
$$p_{x_{i}}=-r_{i}\sin \theta_{i},\qquad p_{y_{i}}=r_{i}\cos \theta_{i},\qquad i=1,2,$$
we parametrize the set $\Lambda_{0}$ with 
$(r_{1},r_{2},\theta_{1},\theta_{2},w)$ such that $r_{1},r_{2}\geq 0, r_1^2+r_2^2=1,\ \theta_1,\theta_{2}\in S^1, \ w\in \R$.
Solving the Hamiltonian system defined by $H$ one can show that the arclength geodesic $\g(t)=(x_{1}(t),x_{2}(t),y_{1}(t),y_{2}(t),z(t))$ associated with an initial covector $p_{0}=(r_{1},r_{2},\theta_{1},\theta_{2},w)\in \Lambda_{0}$ and $|w|\neq 0$, is described by the equations (we restrict  to the case $w>0$, the case $w<0$ is analogous by symmetry)
\begin{align} \label{eq:exp22}
x_i(t)&=\frac{r_i}{\al_i w}(\cos (\al_i w t +\theta_i) - \cos
\theta_i),\notag \\
y_{i}(t)&=\frac{r_i}{\al_i w}(\sin (\al_i w t +\theta_i) - \sin\theta_i), \qquad i=1,2,\\
z(t)&=\frac{1}{2w^2}\left(wt-\sum_{i=1}^{2} \frac{r_i^2}{\al_i} \sin \al_i w t \right). \notag
\end{align}
If $w=0$, geodesics are straight lines contained in the hyperplane $\{z=0\}$
\begin{align} \label{eq:exp234}
x_i(t)&=r_{i}t\cos\theta_{i} \\
y_{i}(t)&=r_{i}t\sin\theta_{i}, \notag\\
z(t)&=0. \notag
\end{align}

From equations \eqref{eq:exp22} one easily shows that the projection of a non-horizontal extremal on every 2-plane $(x_i,y_i)$ is a circle, with period $T_i$, radius $\rho_i$ and center $C_i$ defined by 
\begin{equation} \label{eq:circ}
T_{i}=\frac{2\pi}{\al_i w}, \qquad\rho_i=\frac{r_i}{\al_iw} \qquad C_i=-\frac{r_i}{\al_iw}(\cos \theta_i,\sin \theta_i), \qquad \all i=1, \ldots,\ell.
\end{equation} 
Moreover, generalizing the analogous  property of the 3D Heisenberg group $\mb{H}$, one can recover that the $z$-component of the extremal at time $t$ is the weighted sum (with coefficients $\al_i$) of the areas spanned by the vectors $(x_{i}(t),y_{i}(t))$ in $\R^{2}$  (see also Figure \ref{f:2H}).

We introduce the functions $\rho_{i}(t)=\sqrt{x_{i}(t)^{2}+y_{i}(t)^{2}}$ that satisfy
\begin{align}\label{eq:rho}
\rho_{i}(t)=\frac{2 r_{i}}{\al_{i} |w|}\sin\left(\frac{\al_i w t}{2}\right)= r_{i}\, t \, \sin_{c}\left(\frac{\al_i w t}{2}\right), \quad \text{where}\quad \sin_{c}(x)=\frac{\sin x}{x}.\end{align}
\begin{center}
\begin{figure}[!]
\hspace{1.5cm}\scalebox{1} 
{
\begin{pspicture}(0,-3.2930768)(13.4818945,2.9196186)
\definecolor{color309b}{rgb}{0.8,0.8,0.0}
\definecolor{color309g}{rgb}{0.8,1.0,0.2}
\definecolor{color382b}{rgb}{0.4,0.8,1.0}
\definecolor{color382g}{rgb}{0.2,0.8,1.0}
\rput{-137.50327}(17.75112,8.969973){\psarc[linewidth=0.04,fillstyle=gradient,gradlines=2000,gradbegin=color382g,gradend=color382b,gradmidpoint=1.0,fillcolor=color382b](10.619529,1.0337609){1.5707922}{352.97665}{139.818}
\psline[linewidth=0.04](12.178534,0.8416939)(9.419445,2.0472639)}
\rput{-137.50327}(5.161716,4.233064){\psarc[linewidth=0.04,fillstyle=gradient,gradlines=2000,gradbegin=color309g,gradend=color309g,gradmidpoint=1.0,fillcolor=color309b](3.4038627,1.112976){1.2825377}{6.5215425}{238.08133}
\psline[linewidth=0.04](4.6781015,1.2586424)(2.725766,0.02435868)}
\psline[linewidth=0.03cm,arrowsize=0.05291667cm 2.0,arrowlength=1.4,arrowinset=0.4]{->}(2.54,-2.2553813)(2.56,2.9046187)
\psline[linewidth=0.03cm,arrowsize=0.05291667cm 2.0,arrowlength=1.4,arrowinset=0.4]{->}(0.0,0.14461866)(5.44,0.14461866)
\psline[linewidth=0.03cm,arrowsize=0.05291667cm 2.0,arrowlength=1.4,arrowinset=0.4]{->}(9.32,-2.2753813)(9.34,2.8846188)
\psline[linewidth=0.03cm,arrowsize=0.05291667cm 2.0,arrowlength=1.4,arrowinset=0.4]{->}(6.78,0.124618664)(12.22,0.124618664)
\usefont{T1}{ptm}{m}{n}
\rput(6.321455,-3.0703814){$z(t)=\alpha_1\text{Area}(A_1(t))+\alpha_2 \text{Area}(A_2(t))$}
\usefont{T1}{ptm}{m}{n}
\rput(3.751455,2.6896186){$(x_1(t),y_1(t))$}
\usefont{T1}{ptm}{m}{n}
\rput(12.211455,1.3896186){$(x_2(t),y_2(t))$}
\usefont{T1}{ptm}{m}{n}
\rput(3.801455,0.86961865){$A_1(t)$}
\usefont{T1}{ptm}{m}{n}
\rput(10.681455,-0.21038133){$A_2(t)$}
\end{pspicture} 
}
\caption{Projection of a non-horizontal geodesic: case $0<\alpha_{1}\leq \alpha_{2}$}
\label{f:2H}
\end{figure}
\end{center}

\subsubsection{Non-contact case} In the  case when $\alpha_{2}=0$ the sub-Riemannian structure reduces to the product  $\mb{H}\times\R^{2}$. The level set of the Hamiltonian in this case is again
 $$\Lambda_{0}=\{(p_{x_{1}},p_{y_{1}},p_{x_{2}},p_{y_{2}},w)\,|\, p^{2}_{x_{1}}+p^{2}_{y_{1}}+p^{2}_{x_{2}}+p^{2}_{y_{2}}=1\}\simeq S^{3}\times \R$$
Performing the change of variable
$$p_{x_{i}}=-r_{i}\sin \theta_{i},\qquad p_{y_{i}}=r_{i}\cos \theta_{i},\qquad i=1,2,$$
we parametrize the set $\Lambda_{0}$ with 
$(r_{1},r_{2},\theta_{1},\theta_{2},w)$ such that $r_{1},r_{2}\geq 0, r_1^2+r_2^2=1,\ \theta_1,\theta_{2}\in S^1, \ w\in \R$.
 
Again, solving the Hamiltonian system defined by $H$, one  shows that the arclength geodesic $\g(t)=(x_{1}(t),x_{2}(t),y_{1}(t),y_{2}(t),z(t))$ associated with an initial covector $p_{0}=(r_{1},r_{2},\theta_{1},\theta_{2},w)\in \Lambda_{0}$ and $|w|\neq0$,  are described as follows (we restrict  to the case $w>0$, the case $w<0$ is analogous by symmetry): when $r_{1}\neq 0$
\begin{align} \label{eq:exp23}
x_1(t)&=\frac{r_{1}}{\al_1 w}(\cos (\al_1 w t +\theta_{1}) - \cos
\theta_{1}),\notag \\
y_{1}(t)&=\frac{r_{1}}{\al_1 w}(\sin (\al_1 w t +\theta_{1}) - \sin\theta_{1}),\\
x_2(t)&= r_{2}t \cos \theta_{2} ,\notag \quad
y_{2}(t)=r_{2}t \sin \theta_{2}, \notag\\
z(t)&=\frac{r_{1}^{2}}{2w^2}\left(wt- \frac{1}{\al_1} \sin (\al_1 w t) \right). \notag
\end{align}
If $r_{1}=0$ (hence $r_{2}=1$) we have
\begin{align} \label{eq:exp234}
x_1(t)&=0,\quad
y_{1}(t)=0, \\
x_2(t)&= t\cos\theta_{2} , \quad
y_{2}(t)=t\sin \theta_{2}, \notag\\
z(t)&=0. \notag
\end{align}
One can show that the trajectories corresponding to the case to $r_{1}=0$ are also abnormal. 
Finally, if $w=0$, we have
\begin{align} \label{eq:exp234}
x_1(t)&=r_{1}t\cos\theta_{1},\quad
y_{1}(t)=r_{1}t\sin\theta_{1}, \\
x_2(t)&= r_{2}t\cos\theta_{2}, \quad
y_{2}(t)=r_{2}t\sin\theta_{2}, \notag\\
z(t)&=0. \notag
\end{align}

Let us denote the exponential map starting from the origin as the map
$$\EXP_{0}: \Lambda_{0}\times \R^{+} \to M, \qquad \EXP_0(p_0,t)=\g(t),$$
where $\g(t)=(x_{1}(t),x_{2}(t),y_{1}(t),y_{2}(t),z(t))$ is the arclength geodesic associated with $p_{0}$.
 
 The following lemma is proved in \cite{corank1}, see also \cite{bealsgaveaugreiner}.
\bl An arclength geodesic $\g(t)=\EXP_0(p_0,t)$ associated to $p_{0}=(r_{1},r_{2},\theta_{1},\theta_{2},w)\in \Lambda_{0}$ is optimal up to its cut time
\begin{equation}\label{eq:bicut}
\tcut(\g)=\frac{2\pi}{|w|\max\{\al_1,\al_2\}}.
\end{equation}
with the understanding $\tcut(\g)=+\infty$ when $w=0$. Moreover the cut time coincides with the first conjugate time.
\el

Let us mention the following explicit formula for the distance from the origin to the vertical axis in terms of the parameters $\al_{1},\al_{2}$, namely for every $\zeta=(0,0,0,0,z)$ we have
\begin{equation}\label{eq:dist}
 d(0,\zeta)^2=\frac{4 \pi |z|}{\max\{\al_1,\al_2\}}. 
\end{equation}

\subsection{Proof of Theorem \ref{t-cut}}
By symmetry we will consider the case $w>0$ and the cut locus will be the union of this set and its symmetric with respect to the hyperplane $\{z=0\}$. (Recall the case $w=0$ there is no cut locus along the geodesic.) 

(i). In this case let $t_{*}=2\pi/\al w$ be the cut time, where $\al:=\al_{1}=\al_{2}$. Substituting the cut time into the equations of the geodesic one gets that all the horizontal coordinates vanishes and $z(t_{*})=\pi/\al w^{2}$. From this (i) is immediate. Notice moreover that the final point does not depend on $(r_{1},r_{2},\theta_{1},\theta_{2})\in S^{3}$, so we have a three-parameter family of optimal geodesics joining this point. \\[0.2cm]
(ii). Since $\max\{\al_{1},\al_{2}\}=\al_{1}$, the cut time is $t_{*}=2\pi/\al_{1} w$ and substituting the cut time into the equations we get
\begin{align} \label{eq:exp22cut}
x_{1}(t_{*})&=y_{1}(t_{*})=0,\nn\\
x_2(t_{*})&=\frac{r_2}{\al_2 w}\left(\cos \left(2\pi \frac{\al_2}{\al_{1}}+\theta_2\right) - \cos
\theta_2\right),\notag \\
y_{2}(t_{*})&=\frac{r_2}{\al_2 w}\left(\sin \left(2\pi \frac{\al_2}{\al_{1}} +\theta_2\right) - \sin\theta_2\right), \\
z(t_{*})&=\frac{1}{2w^2}\left(\frac{2\pi}{\al_{1}}-\frac{r_2^2}{\al_2} \sin \left(2\pi \frac{\al_2}{\al_{1}}\right) \right). \notag
\end{align}
Notice that the set of cut points has rotational symmetry. Indeed the function $\rho_{2}$ describing the distance from the origin in the plane $(x_{2},y_{2})$ (see \eqref{eq:rho})
$$\rho_{2}^{2}=x_{2}^{2}+y^{2}_{2}=\frac{4 r^{2}}{\al_{2}^{2} w^{2}}\sin^{2}\left(\pi\frac{\al_{2}}{\al_{1}}\right)$$
does not depend on $\theta_{2}\in S^{1}$. In other words the cut locus  is a subset of the three-dimensional space $\{(x_{2},y_{2},z)\}$ which is symmetric with respect to the origin. This set is obtained by rotating with respect to the $z$ axis the image of the map in the half-space $\{(\rho,z),\rho\geq 0\}$ (for simplicity we denote in the following by $\rho,r$ respectively $\rho_{2}$ and $r_{2}$)
\begin{equation} \label{eq:exp22cut2}
(r,w)\mapsto
\begin{cases}
\displaystyle{\rho= \frac{2 r}{\al_{2} w}\sin\left(\pi\frac{\al_{2}}{\al_{1}}\right)} \\[0.3cm]
\displaystyle{z=\frac{1}{2w^2}\left(\frac{2\pi}{\al_{1}}-\frac{r^2}{\al_2} \sin \left(2\pi \frac{\al_2}{\al_{1}}\right)\right)}  
\end{cases}, \qquad r\in[0,1], \,w>0.
\end{equation}
From the equations above it is easy to get the relation
$$z=\Psi(\al_{1},\al_{2},r)\rho^{2}$$
where $\Psi(\al_{1},\al_{2},r)$ is a constant that depends only on $r$ (and the parameters $\al_{1},\al_{2}$) expressed by
$$\Psi(\al_{1},\al_{2},r)=\frac{\al_{2}^{2}}{8 r^{2}}\sin^{-2}\left(\pi \frac{\al_{1}}{\al_{2}}\right)\left(\frac{2\pi}{\al_{1}}-\frac{r^{2}}{\al_{2}}\sin\left(2\pi\frac{\al_{2}}{\al_{1}}\right)\right) . $$
In other words, for a fixed value of $r$, the point belongs to a parabola and \eqref{eq:exp22cut2} guarantees there exists a unique $w$ associated to each point on this parabola.
Moreover we have that
$$\lim_{r\to 0^{+}}\Psi(\al_{1},\al_{2},r)=+\infty,\qquad \lim_{r\to 1^{-}}\Psi(\al_{1},\al_{2},r)=:K(\al_{1},\al_{2}),$$
and
$$\frac{\partial}{\partial r}\Psi(\al_{1},\al_{2},r)=-\frac{\pi\al_{2}^{2}\csc^{2}(\pi\frac{\al_{1}}{\al_{2}})}{8\al_{1}r^{3}}<0, \qquad r\in[0,1],$$
that ensure that the function \eqref{eq:exp22cut2} is injective on its image, which completes the proof of (ii).
 
(iii). This case is reduced to $\R^{2}$ times the cut locus of $\mb{H}$, that is the three dimensional space $\{(x_{2},y_{2},z), z\neq 0\}$.

\brem 
Notice that the structure of the cut locus is continuous with respect to the parameters $\al_{1},\al_{2}$ in the following sense
$$\lim_{\al_{2}\to 0^{+}}\Psi(\al_{1},\al_{2},r)=0,\qquad \lim_{\al_{2}\to \al_{1}^{-}}\Psi(\al_{1},\al_{2},r)=+\infty.$$
This implies in fact that the cases (i) and (iii) can be recovered as the limit cases of (ii). 
\erem

\subsection{Proof of Theorem \ref{t-heat}}
Point (i) and (ii) of Theorem \ref{t-heat} follows from results of Ben Arous \cite{benarouscut,benarousdiag}.

Point (iii) of Theorem \ref{t-heat} can be obtained by applying Theorem \ref{t-main} to the explicit structure of optimal geodesics described in the previous section. We have two different cases: 

(iii.a).  $0< \al_{2}=\al_{1}=:\al$. In this case we proved that every point in the cut locus is reached by a three-parameter family of geodesics. More precisely, following the notation of Theorem \ref{t-main}, for every point $q\in \Cut_{0}\setminus\{0\}$, we have $\mathcal{O}=S^{3}$. Moreover it is easy to see that for every $\bar p\in \mathcal O$ one has $D_{(\bar p,d(0,q))}\Exp_{0}$ has rank 2. Indeed the differential of the exponential map is never degenerate with respect to $t$ and, for $t=d(0,q)$ fixed, one has (notice that $\bar p\in \mathcal O$ implies $w\neq 0$)
$$\Exp_{0}(\bar p,t)=(0,0,0,0,\pi/ (\al w^{2})).$$
Hence, applying Theorem \ref{t-main} with $n=5$ and $r=3$, we have
$$p_{t}(0,q)=\frac{C+O(t)}{t^{4}}e^{-\frac{d^{2}(0,q)}{4t}}$$

(iii.b).  $0\leq \al_{2}<\al_{1}$. Let us first consider the subcase $\al_{2}>0$. Then every point in the cut locus is reached by a one-parameter family of geodesics. More precisely, following the notation of Theorem \ref{t-main}, for every point $q\in \Cut_{0}\setminus\{0\}$, we have $\mathcal{O}=S^{1}$. Moreover it is easy to see that for every $\bar p\in \mathcal O$ one has $D_{(\bar p,d(0,q))}\Exp_{0}$ has rank 4. 
Hence, applying Theorem \ref{t-main} with $n=5$ and $r=1$, we have
$$p_{t}(0,q)=\frac{C+O(t)}{t^{3}}e^{-\frac{d^{2}(0,q)}{4t}}$$
When $\al_{2}=0$ the results of Section \ref{s:heat} cannot be applied due to the presence of abnormal minimizers. Nevertheless since the struture is the product of the Heisenberg group and $\R^{2}$ one easily gets that the heat kernel is the product of the heat kernel of the Heisenberg group and the heat kernel of $\R^{2}$. From the explicit expression one gets the result.
%
%

\section{What one can get from the Greiner-Gaveau-Beals results}\label{s:final}
In this section we recall the expression of the heat kernel of the intrinsic sub-Laplacian associated with a  2-step nilpotent structure, that has been found in \cite{bealsgaveaugreiner}. 
Then we rewrite it to have a convenient expression on the \virg{vertical subspace}.

Consider on $\R^n$ a 2-step nilpotent structure of rank $k<n$, where $X_{1},\ldots,X_{k}$ is an orthonormal frame. Once a smooth complement $\mc{V}$ for the distribution is chosen (i.e.\ $T_{q}\R^n=\distr_{q}\oplus \mc{V}_{q}$, for all $q\in \R^n$) we can complete an orthonormal frame to a global one $X_{1},\ldots,X_{k},Y_{1},\ldots,Y_{m}$, where $m=n-k$ and $\mc{V}_{q}=\tx{span}_{q}\{Y_{1},\ldots,Y_{m}\}$. Since the structure is nilpotent, we can assume that the only nontrivial commutation relations are
\begin{equation} \label{eq:B}
[X_{i},X_{j}]=\sum_{h=1}^{m}b_{ij}^{h}Y_{h},
\end{equation}
where $B_{1},\ldots,B_{m}$ defined by $B_{h}=(b_{ij}^{h})$ are skew-symmetric matrices (see \cite{corank2} for the role of these matrices in the exponential map).

Due to the group structure, the intrinsic sub-Laplacian takes the form of sum of squares $\lapl=\sum_{i=1}^{k}X_{i}^{2}$ (see Remark \ref{r:group}). The group structure also implies that the heat kernel is invariant with respect to the group operation,
hence it is enough to consider the heat kernel $p_{t}(0,q)$ starting from the identity of the group, which we also denote $p_{t}(q)$. The heat kernel is written as follows (see again \cite{bealsgaveaugreiner,lanconellibook})
\begin{equation}\label{eq:gaveau}\nn
\ \ p_{t}(q)=\frac{2}{(4\pi t)^{Q/2}} \int_{\R^{m}} V(B(\tau)) \exp \left(-\frac{ W(B(\tau)) x\cdot x}{4t}\right) \cos \left(\frac{z \cdot \tau}{t}\right) d\tau,
\end{equation}
where  $q=(x,z)$, $x\in \R^k, z\in \R^m$, and $B(\tau):=\sum_{i=1}^{m} \tau_{i}B_{i}$. Moreover  $V: \R^{n\times n} \to \C$ and $W: \R^{n\times n} \to  \R^{n\times n}$ are the matrix functions defined by
$$V(A)=\sqrt{\det \left(\frac{A}{\sin A}\right)}, \qquad W(A)= \frac{A}{\tan A}.$$
Here $Q$ is the Hausdorff dimension of the sub-Riemannian structure.

Notice that \eqref{eq:gaveau} differs by some constant factors from the formulas contained in  \cite{bealsgaveaugreiner} since there the heat kernel is the solution of the equation $\partial_t u=\frac12 \lapl u$. 
\brem
Assume that the real skew-symmetric matrix $B(\tau)$ is diagonalizable and denote by $\pm i \lam_{j}(\tau)$, for $j= 1,\ldots \ell$, its non-zero eigenvalues. Then we have the formula for the expansion on the \virg{vertical subspace} (i.e.\ where $x=0$)
\begin{equation} \label{eq:z}
p_{t}((0,z))=\frac{2}{(4\pi t)^{Q/2}} \int_{\R^{m}} \prod_{j=1}^{\ell} \frac{\lam_{j}(\tau)}{\sinh \lam_{j}(\tau)}\cos \left(\frac{z \cdot \tau}{t}\right) d\tau.
\end{equation}
\erem

For the bi-Heisenberg case, in which $m=1$, we have that $B(\tau)=\tau B$ and $\text{eig}(B(\tau))=\{\pm\, i\al_1\tau,\pm\, i\al_2 \tau\}$, and from \eqref{eq:gaveau} one gets
\begin{align}
p_{t}(q)=
\frac{2}{(4\pi t)^3}&\int_{-\infty}^{\infty} \prod_{i=1}^{2}\frac{\al_i\tau}{\sinh(\al_i\tau)}\exp\left(-\sum_{i=1}^{2}
\frac{x_i^2+y_i^2}{4 t} \frac{\al_i\tau}{\tanh(\al_i\tau)}
\right)\cos\left(\frac{z\tau}{t}\right) d\tau,\nn
\end{align}
where $p_{t}(q):=p_{t}(0,q)$ and $q=(x_{1},x_{2},y_{1},y_{2},z)$.

From this formula one can get the expansion of the heat kernel on the $z$-axis (that always lies in the cut locus by Theorem  \ref{t-cut}). We rewrite the expansion when $x_{i}=y_{i}=0$ for $i=1,2$, i.e.\ at a point $ \zeta=(0,0,0,0,z)$ (see also \eqref{eq:z}), as
\begin{equation} \label{eq:a1ua2}
p_{t}(\zeta)=\frac{2}{(4\pi t)^3}\int_{-\infty}^{\infty}\frac{\al_1\tau}{\sinh(\al_1\tau)}\frac{\al_2\tau}{\sinh(\al_2\tau)}
\cos\left(\frac{z\tau}{t}\right) d\tau. 
\end{equation}
We show now that the different behavior of the asymptotics that we proved in the previous section can be recovered in the two cases when the parameters are equal or not. 
\bi
\iii[(a)]\emph{Case $\al_1=\al_2$.} It is enough to consider the case when $\alpha_{1}=\alpha_{2}=1$. 

In this case the integral \eqref{eq:a1ua2} is explicitly computed by the formula
\begin{align}\label{eq:expl1}
p_{t}(\zeta)&=\frac{2}{(4\pi t)^3}\int_{-\infty}^{\infty}\frac{\tau^{2}}{\sinh^{2}\tau}
\cos\left(\frac{z\tau}{t}\right) d\tau=\frac{\pi  z \coth \left(\frac{\pi  z}{2 t}\right)-2 t}{32 \pi 
   t^4 \left(\cosh \left(\frac{\pi  z}{t}\right)-1\right)}.
   \end{align}
Let us consider now the point $\zeta$ corresponding to $z=1$. From \eqref{eq:dist} and \eqref{eq:expl1} one gets    
\begin{equation} \label{eq:as2}
   p_t(\zeta)=\frac{1}{t^4} \exp\left(-\frac{d^2(0,\zeta)}{4t}\rp
\smooth(t),
   \end{equation}
   with $\smooth$ a smooth function of $t$ such that $\smooth(0)=1/16$.

\iii[(b)] \emph{Case $\al_1\neq \al_2$.}  To simplify the expression of the heat kernel we consider the particular case when $\al_1=1, \al_2=1/2$. Also in this case we have the  expression of the heat kernel on the vertical axis $\zeta=(0,0,0,0,z)$ (again $z\neq0$)
\begin{equation} \label{eq:expl2}
\qquad  p_t(\zeta)=\frac{ \left(-8 \pi  z \sinh
   \left(\frac{\pi  z}{t}\right)+8 t \cosh \left(\frac{\pi 
   z}{t}\right)+\pi  t \left(\cosh \left(\frac{2 \pi 
   z}{t}\right)-3\right)\right)}{128 \pi  t^4  \text{cosh}^3\left(\frac{\pi  z}{t}\right)}.
\end{equation}
If we consider the expansion at the point  $\zeta$ corresponding to $z=1$, using \eqref{eq:dist} and \eqref{eq:expl2}  we find 
\begin{equation} \label{eq:as1} 
p_t(\zeta)=\frac{1}{t^3}  \exp \lp-\frac{d^2(0,\zeta)}{4t}\rp
   \smooth(t),
   \end{equation}
   with $\smooth(0)=1/32$.
\ei
By homogeneity these expansion hold on any point of the $z$-axis.
\brem 
As a byproduct of formula \eqref{eq:as2}, we get that, in the symmetric case the maximal degeneration, at least in terms of the factor in front of the exponential, is not obtained on the diagonal. This is in contrast to what happens in the Heisenberg case.
\erem


{\bf Acknowledgements.}
This research has been partially supported by the European Research Council, ERC StG 2009 “GeCoMethods”, contract n. 239748, by the iCODE institute (research project of the Idex Paris-Saclay), by the ANR project SRGI  ``Sub-Riemannian Geometry and Interactions'', contract number ANR-15-CE40-0018, and by the National Security Agency under Grant Number H98230-15-1-0171. This research benefited also from the support of the “FMJH Program Gaspard Monge in optimization and operation research” and from the support to this program from EDF.

{\small
\bibliography{2Heis-biblio}
\bibliographystyle{abbrv}
}

\end{document}